\documentclass{article}

\usepackage{PRIMEarxiv}

\usepackage[utf8]{inputenc} % allow utf-8 input
\usepackage[T1]{fontenc}    % use 8-bit T1 fonts
\usepackage{hyperref}       % hyperlinks
\usepackage{url}            % simple URL typesetting
\usepackage{booktabs}       % professional-quality tables
\usepackage{amsfonts}       % blackboard math symbols
\usepackage{nicefrac}       % compact symbols for 1/2, etc.
\usepackage{microtype}      % microtypography
\usepackage{lipsum}
\usepackage{fancyhdr}       % header
\usepackage{graphicx}       % graphics
\usepackage{adjustbox}
\graphicspath{{media/}}     % organize your images and other figures under media/ folder

\title{Power Availability of PV plus Thermal Batteries in real-world electric power grids}

\author{Odin Foldvik Eikeland$^{1,2}$, Colin C. Kelsall$^{1}$, Kyle Buznitsky$^{1}$, Shomik Verma$^{1}$, Filippo Maria Bianchi$^{3}$, \\
\textbf{Matteo Chiesa$^{1,2,4}$\thanks{Corresponding author: mchiesa@mit.edu}, \& Asegun Henry$^{1}$\thanks{Corresponding author: ase@mit.edu}} \\
\\
$^1$Department of Mechanical Engineering, Massachusetts Institute of Technology, Cambridge, MA, USA.\\
$^2$Department of Physics and Technology, UiT – the Arctic University of Norway, 9037 Tromsø, Norway.\\
$^3$Department of Mathematics and Statistics, UiT – the Arctic University of Norway Technology, 9037 Tromsø, Norway.\\
$^4$Laboratory for Energy and NanoScience (LENS), Khalifa University of Science and Technology, Masdar.\\
}

\begin{document}

\maketitle
%\footnote[1]{Department of Mechanical Engineering, Massachusetts Institute of Technology, Cambridge, MA, USA.}
%\footnote[2]{Department of Physics and Technology, UiT – the Arctic University of Norway, 9037 Tromsø, Norway.}
%\footnote[3]{Department of Mathematics and Statistics, UiT – the Arctic University of Norway Technology, 9037 Tromsø, Norway.}
%\footnote[4]{Laboratory for Energy and NanoScience (LENS), Khalifa University of Science and Technology, Masdar.}

\begin{abstract}
As variable renewable energy sources comprise a growing share of total electricity generation, energy storage technologies are becoming increasingly critical for balancing energy generation and demand.

In this study, we modeled an existing thermal energy storage unit with estimated capital costs that are sufficiently low to enable large-scale deployment in the electric power system. Our analysis
emphasizes the value of using such units to cost-effectively improve renewable energy dispatchability.
This study modeled an existing real-world grid rather than simulating hypothetical future electric power systems. The storage unit coupled with a photovoltaic (PV) system was modeled with different storage capacities, whereas each storage unit size had various discharge capacities. 

The modeling was performed under a baseline case with no emission constraints and under hypothetical scenarios in which CO$_2$ emissions were reduced. The results show that power availability increases with increasing storage size and vastly increases in the hypothetical CO$_2$ reduction scenarios, as the storage unit is utilized differently. When CO$_2$ emissions are reduced, the power system must be less dependent on fossil fuel technologies that currently serve the grid, and thus rely more on the power that is served from the PV + storage unit. 

The proposed approach can provide increased knowledge to power system planners regarding how adding PV + storage systems to existing grids can contribute to the efficient stepwise decarbonization of electric power systems.

\end{abstract}

% keywords can be removed
\keywords{Thermal Energy Grid Storage \and Power System Modeling \and Decarbonizing Electric Grids}

\newpage
\section{Introduction}
The use of variable renewable energy (VRE) resources, such as wind power and solar photovoltaics (PV), is expanding rapidly as a share of total power generation and is critical to the decarbonization of electrical power systems \cite{Eriksen2022,gielen2019role,apostoleris2019utility}. The weather-dependent intermittency of VRE sources complicates the planning and management of power systems as the electric power generation can no longer be directly modulated to match the electricity demand. Energy storage will therefore be an increasingly critical component of future energy systems with high penetrations of VRE sources. Energy storage can charge excess electricity in periods with high generation and low demand, and then discharge the electricity in periods with low generation from the VRE sources to match the load in periods with high electricity demand \cite{braff2016value,ziegler2019storage}.

The need for inexpensive storage over periods with different lengths, from seconds to days and even seasonal storage, has accelerated in accordance with the increasing share of VRE technologies in electric power systems \cite{shaner2018geophysical,sepulveda2021design}. Pumped hydropower storage (PHS) is an established technique for large-scale energy storage but can be used only in certain geographical areas. Lithium-ion batteries have been the state-of-the-art technology for short-term storage. However, capital costs between US\$80 and US\$100 kWh$^{-1}$ make them unaffordable for the multi-day storage objectives required to completely decarbonize the grid \cite{braff2016value,ziegler2019storage,mallapragada2020long}.  Concentrated solar power with thermal energy storage (CSP-TES) has been seen as a promising option, but major projects around the world have been plagued by delays, cost overruns and mechanical issues, and interest has waned in recent years \cite{xu2016prospects,kennedy2022role}.  Studies suggest that achieving cost-efficient multi-day storage requires a capital cost reduction to US\$3-30 kWh$^{-1}$ \cite{ziegler2019storage,albertus2020long}. Resolving this issue could enable more rapid decarbonization of the power system, resulting in a 25\% reduction in global GHG emissions \cite{IPCC2014,EPA2011}. Therefore, one of the most significant issues that needs to be resolved to achieve the GHG emission reduction targets is the energy storage challenge.

A storage concept based on Thermal Energy Storage (TES) has shown promising potential to achieve sufficiently low capital cost in the multi-day storage regime. TES stores the electricity as heat rather than electrochemically, and then converts it back to electricity when needed \cite{henry2020five}. The Thermal Energy Grid Storage (TEGS) concept, detailed in \cite{amy2019thermal}, stores electricity as sensible heat in graphite storage blocks and uses thermophotovoltaics (TPV) to convert heat back to electricity on demand \cite{amy2019thermal,kelsall2021technoeconomic}. While the conversion of heat to electricity results in significant efficiency penalties, storing energy as heat instead of electrochemically can be vastly cheaper, and thus the round-trip efficiency (RTE) penalty compared to electrochemical batteries ($\sim$ 90\%) can potentially be a worthwhile tradeoff \cite{amy2019thermal}. To maximize the conversion efficiency from heat to electricity, the heat is stored at extremely high temperatures (~2400 °C). In a recent work by \cite{lapotin2022thermophotovoltaic}, the authors demonstrated a world-record high conversion efficiency of 41\% using TPV, and reported a projected conversion efficiency of 50\% in the future. As such this technology can achieve a projected cost below US\$ 20 kWh$^{-1}$ at gigawatt scales.

In addition to the projected low cost, a unique property of TEGS compared to Li-ion battery technology is the fact that, since energy is stored as heat in graphite blocks and thereafter converted to electricity using TPV, it enables the possibility of fully decoupling the charge and discharge capacities of the storage unit. 
This allows the TEGS to charge (i.e., store heat) at a much higher capacity than that required for discharging. The benefit of such a mechanism is that a large amount of energy can be charged in a short amount of time when generation surpluses exist and discharged over a longer period to cover the electricity load in periods where demand exceeds supply.

The body of existing literature counts several studies that have employed different approaches to evaluate the value of using storage to increase the dispatchability of VRE sources, and the different studies have highlighted the storage requirements (capital cost and storage duration) to enable the full decarbonization of the power system. Tab. 1 shows an overview of some relevant studies, where the key findings in each work are highlighted. 

Despite the vast number of previous studies on modeling the value of energy storage in emerging power systems, there is a lack of research that analyzes specific real-world grids. Previous studies have largely focused on modeling hypothetical future electric power systems starting from scratch (i.e., “greenfield” models) \cite{sepulveda2021design,schill2018long,de2016value,heuberger2017systems,sepulveda2018role,frew2016flexibility,liu2019role}. However, such studies can, in many cases, lead to vague results as this requires a complete change of the current electricity mix. This is in many cases challenging to implement due to policy considerations. In addition, by modeling greenfield cases, infrastructure that already exists is neglected.

Some studies have modeled existing grids (i.e., brownfield studies), but lack in modeling the value of using specific storage technologies \cite{ziegler2019storage,mallapragada2020long,jafari2020power,schleifer2022exploring}. These studies model the storage requirements in general, whereas all studies show that the studies show that the capital cost must be below US \$20 kWh$^{-1}$. No previous work has modeled the potential of using an emerging storage concept based on TES that already exists on a lab scale \cite{amy2019thermal,kelsall2021technoeconomic}.

In this study, we propose a framework for addressing the value of using TEGS that has sufficiently low capital cost to be economically used in an electric grid. Using a Capacity Expansion Model (CEM) \cite{jenkins2017enhanced}, we consider a hypothetical PV + TEGS system that is interconnected to an existing real-world grid in the Northeastern US. The TEGS unit charges excess electricity from PV during periods of surplus generation. When the grid demands electricity and the PV plant cannot deliver sufficient power due to a lack of solar availability, the stored energy is discharged. Different storage sizes with varying discharge capacities connected to the PV plant are modeled to optimize energy availability. To investigate how emission constraints affect the energy availability of PV + storage systems, a hypothetical future scenario is modeled for the existing power system where CO$_2$ emissions are reduced.

The main contributions of this study are: Rather than modeling the power system as a greenfield case study, we analyze an abstract representation of an existing grid, i.e., a “brownfield” model, and address how adding a PV + storage system can contribute to decarbonizing the grid. Instead of modeling general requirements of storage to enable the full decarbonization of the power system, we here model a TES unit that currently exists at lab-scale and has promising cost projections that are well-documented in the literature \cite{amy2019thermal}. This study can provide increased knowledge to power system planners on how coupling emerging storage technologies and PV systems to existing grids can contribute to stepwise decarbonizing of the grid in a more short-term horizon.

The remainder of this paper is organized as follows. Section 2 presents the share of technologies in the modeled existing electric power system, and the methodology for modeling the grid in which the hypothetical PV + storage system is added to interact with the existing power system. Section 3 presents the metrics for computing the power availability of electric power systems. In Section 4, the results are presented in terms of power availability under different modeling scenarios. Finally, the conclusions are presented in Section 5. 

\begin{table}[]
\centering
\footnotesize
\setlength\tabcolsep{.9em} %horizontal padding
\caption{Overview of relevant work addressing the value of energy storage}
\label{tab:allworks}
\begin{tabular}{llllll}
\multicolumn{1}{l|}{\textbf{Ref.}} & \multicolumn{1}{l|}{\textbf{Year}} & \textbf{Key findings}                                                                                                                                                                                              \\ \hline
\multicolumn{1}{l|}{\cite{de2016value}}              & \multicolumn{1}{l|}{2016}          & Large-scale deployment of available battery technologies requires cost reductions.                                                                                                                                  \\
\multicolumn{1}{l|}{\cite{frew2016flexibility}}              & \multicolumn{1}{l|}{2016}          & Pathways to fully renewable systems are feasible with high cost and overgeneration.                                                                                                                                \\
\multicolumn{1}{l|}{\cite{braff2016value}}              & \multicolumn{1}{l|}{2016}          & Cost reduction for storage technologies is required to reach widespread profitability.                                                                                                                              \\
\multicolumn{1}{l|}{\cite{heuberger2017systems}}              & \multicolumn{1}{l|}{2017}          & The availability of how low-carbon technologies impact the optimal capacity mix and generation patterns were demonstrated.                                                                                          \\
\multicolumn{1}{l|}{\cite{shaner2018geophysical}}              & \multicolumn{1}{l|}{2018}          & To reliably meet 100\% of total annual electricity demand, weeks of energy storage are required to support with electricity.                                                                                        \\
\multicolumn{1}{l|}{\cite{schill2018long}}              & \multicolumn{1}{l|}{2018}          &  \begin{tabular}[c]{@{}l@{}}The role of energy storage units in power systems with high shares of VRE was analyzed.\\ The importance of storage increases with the increasing share of renewable-based power technologies.       \end{tabular}

\\
\multicolumn{1}{l|}{\cite{sepulveda2018role}}              & \multicolumn{1}{l|}{2018}          & \begin{tabular}[c]{@{}l@{}}Firm low-carbon resources consistently lower decarbonized electricity system costs,\\and the availability of firm low-carbon resources reduces costs 10\%–62\% in zero-CO$_2$ cases.           \end{tabular}
\\
\multicolumn{1}{l|}{\cite{liu2019role}}              & \multicolumn{1}{l|}{2019}          & 
\begin{tabular}[c]{@{}l@{}}The benefits of hydropower and storage units were analyzed. \\Three decarbonized power systems with distinct grid expansion strategies were compared. \\ Cutting transmission volume does not increase the total costs.           \end{tabular}
\\

\multicolumn{1}{l|}{\cite{denholm2019timescales}}              & \multicolumn{1}{l|}{2019}          & Curtailment of renewable energy generation can be avoided using energy storage.                                                                                                                                    \\
\multicolumn{1}{l|}{\cite{ziegler2019storage}}              & \multicolumn{1}{l|}{2019}          & Energy storage cost below \$20 kWh$^{-1}$ can enable cost-competitive baseload power.                          \\                                                                                         
\multicolumn{1}{l|}{\cite{guerra2020value}}              & \multicolumn{1}{l|}{2020}          & Hydrogen storage with up to 2 weeks of discharge duration is expected to be cost-effective in future power systems.                                                                                                
\\
\multicolumn{1}{l|}{\cite{schulthoff2021role}}              & \multicolumn{1}{l|}{2020}          & Hydrogen storage enable for sector coupling in real-world power systems.                                                                                                                                           \\
\multicolumn{1}{l|}{\cite{bompard2020electricity}}              & \multicolumn{1}{l|}{2020}          & Electricity triangle assures a consistent framework for the energy transition.                                                                                                                                     \\
\multicolumn{1}{l|}{\cite{mallapragada2020long}}              & \multicolumn{1}{l|}{2020}          & Current Li-ion capital cost exceeds storage value in many instances.                                                                                                                                                \\
\multicolumn{1}{l|}{\cite{jafari2020power}}              & \multicolumn{1}{l|}{2020}          & Decarbonization is less expensive with Energy Storage Systems, given sufficient low-cost assumptions.                                                                                                               \\
\multicolumn{1}{l|}{\cite{sepulveda2021design}}              & \multicolumn{1}{l|}{2021}          & Energy capacity costs must be $\leq$ US \$ 20/kWh to reduce the electricity price by $\geq$ 10\%.                                                                                                                                \\
\multicolumn{1}{l|}{\cite{cole2021quantifying}}              & \multicolumn{1}{l|}{2021}          & Power systems with 100\% RE is possible using existing technologies.                                                                                                                                               \\
\multicolumn{1}{l|}{\cite{denholm2021challenges}}              & \multicolumn{1}{l|}{2021}          & There is a need for analytic tool development to model how to achieve a power system that are 100\% decarbonized.                                                                                                   \\
\multicolumn{1}{l|}{\cite{baik2021different}}              & \multicolumn{1}{l|}{2021}          & Clean firm resources are cost-effective in decarbonizing the grid.                                                                                                                                              \\
\multicolumn{1}{l|}{\cite{gandhi2022catching}}              & \multicolumn{1}{l|}{2022}          & Green hydrogen cost between \$0.79/kg and \$1.94/kg in 2030 can be achieved.                                                                                                                                                     \\
\multicolumn{1}{l|}{\cite{schleifer2022exploring}}              & \multicolumn{1}{l|}{2022}          & The demand for power capacity will drive future adoption of higher battery power capacity.                                                                                                                          \\
   
\end{tabular}
\end{table}

\section{Methods}
\subsection{Modeled electric power system}
\begin{figure}[!ht]
    \centering
    \includegraphics[scale=1]{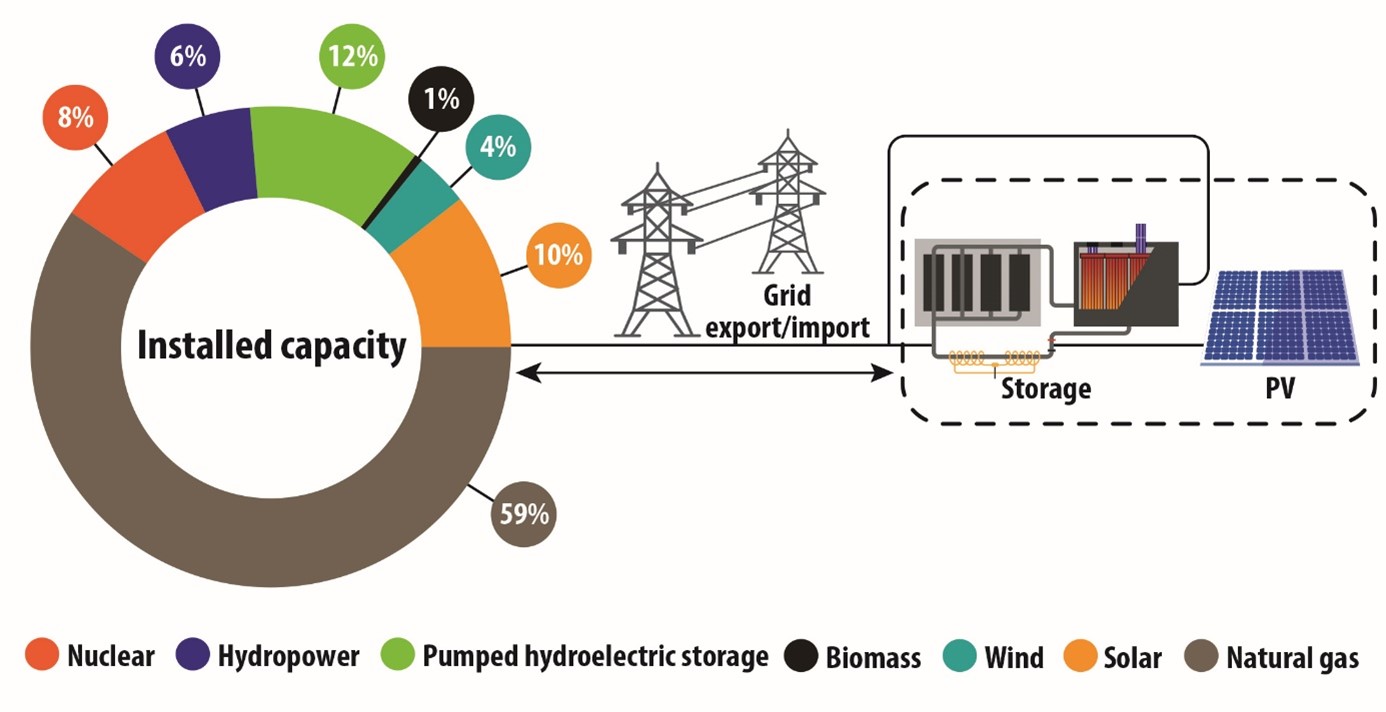}
    \caption{The modeled New England grid zone is interconnected with a hypothetical PV + storage system. The existing power system is dominated by electricity generation from the NG. Solar and wind power represent a small share of electric power systems. To charge the storage unit, electricity (from the grid or PV) is converted to heat using resistance heaters. Heat is stored in graphite blocks at high temperatures (illustrated as gray blocks). To discharge the system, heat is transferred from the graphite blocks using liquid tin to the TPV cells (illustrated as yellow/red pipes to the left of the graphite blocs), where the transported heat is converted to electricity.}
    \label{fig:Gridoverview}
\end{figure}

In this study, an idealized single node representing the electric grid region in the New England power system in North America is modeled. This system considers one grid zone that represents a simplified power system topology of the states of Massachusetts, New Hampshire, and Rhode Island. The electricity demand, capital cost and performance data for the different generation technologies in these regions were collected from the NREL annual technology baseline (ATB) and EIA. The Github library PowerGenome \footnote[1]{The PowerGenome Github library collects source data from EIA, NREL, and EPA and formats the input files for the CEM model. The GitHub library could with associated documentation could be found here \url{https://github.com/PowerGenome/PowerGenome}} was used to collect the input data and shape them to the required format for the CEM. The weather year for modeling VRE availability was 2020. Fig.~\ref{fig:Gridoverview} shows the share of different technologies in the current electric power system. Natural Gas (NG) is the dominant power supply technology, accounting for 59\% of the installed capacity. In the existing power grid, VRE sources such as wind and solar PV represent a smaller share (14\%) of the overall electricity generation mix. 

The hypothetical PV + storage power system is connected to the existing power system through a transmission grid network. Fig.~\ref{fig:Gridoverview} shows how the hypothetical system is connected to the existing grid, where the combined system can fully participate in the power system by exchanging electricity on demand. In this study, we chose to model a storage unit in conjunction with a PV plant, which is believed to be the most dominant source of electricity generation in the future power market \cite{Eriksen2022,gielen2019role,creutzig2017underestimated}. In addition, the normal profile of the daily generation from PV plants is believed to be a good match with storage technologies, as it can store electricity when the PV plant power generates a large amount of electricity during mid-day (and the demand is often low during mid-day) and discharge the stored electricity when the sun is set (during early morning and afternoon/evening). 
In this study, a PV plants with installed peak power capacity of 100 MW and 1 GW were analyzed (see the Supplementary material for the GW scale modeling), while the storage unit was modeled with different energy storage capacities. This system is small compared to the rest of the electric power system and is considered a price taker in this grid. 

\subsection{Capacity Expansion Model (CEM) configuration and storage modeling}
\subsubsection{Capacity Expansion Model (CEM)}
The analysis utilizes GenX \cite{jenkins2017enhanced}, an electric power system CEM that evaluates the cost-optimal electricity mix of the generation, storage, and transmission infrastructure. The cost-optimization is subject to several constraints, such as operational (electricity demand and generation) and policy (CO$_2$ emission) constraints. The dominant constraint used in this study was the maximum limit on the allowed CO$_2$ emissions in the grid. Here, the power system was modeled with a baseline scenario without any CO$_2$ constraints (i.e., the model finds the cost-optimized electricity mix regardless of CO$_2$ emissions) and with a scenario where the CO$_2$ emissions are reduced by 50\%. Constraining the maximum allowed CO$_2$ emissions in the grid will change the dynamics of how the power system is operated as the grid cannot use fossil-fuel-based technologies such as NG at the same scale anymore. To fully capture high-resolution temporal dependencies in the grid, we modeled the grid operation for each hour of the year. All scenarios were evaluated with the Gurobi optimization solver \cite{lubin2015computing} using 16 cores with 128 GB RAM.

\subsubsection{Storage modeling}
Three different TEGS sizes coupled with the 100 MW PV plant were modeled. The modeled energy storage capacities were as follows: 1) 400 MWh, 2) 600 MWh, and 3) 800 MWh. For each storage size, the charging capacity (i.e., the amount of power that can be charged within one hour) is 100 MW, and the discharging capacity is in the range of [5, 100] MW. These storage configurations were modeled under the baseline case (with no CO$_2$ reduction constraint) and 50\% CO$_2$ reduction scenarios. In total, 66 scenarios were modeled to address PV + TEGS energy availability with different storage sizes, discharge capacities, and CO$_2$ constraints. Because future grid scenarios are modeled in this study, the TEGS system is assumed to have a 50\% RTE.

\section{Experimental evaluations}
\subsection{Computing the availability of electric power systems}
\begin{figure}[!ht]
    \centering
    \includegraphics[scale=1]{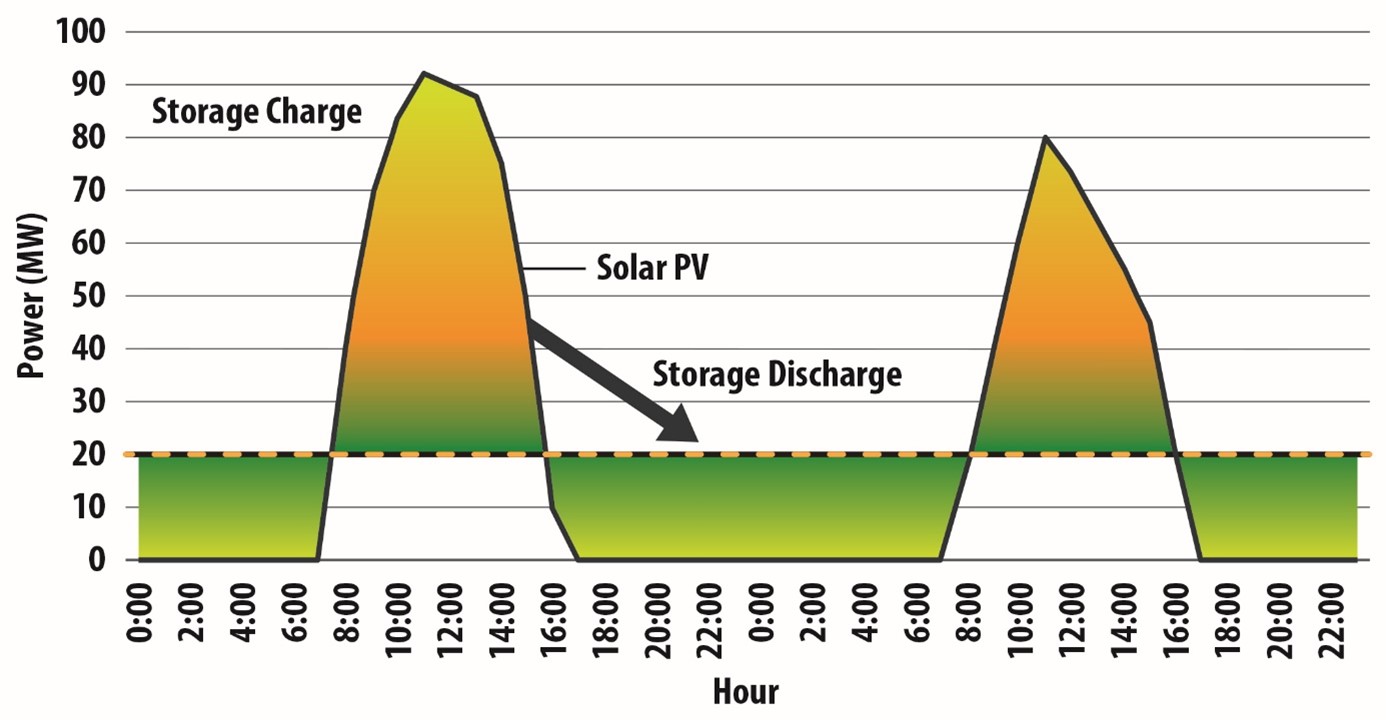}
    \caption{Example illustration where storage is used to provide a constant baseload power to the grid by charging when there is excess electricity and discharging when the PV plant does not generate electricity. The figure also illustrates the large variability in the daily generation profile of solar PV systems. }
    \label{fig:TEGS_illustration}
\end{figure}

In this study, we are interested in assessing the amount of time our hypothetical PV + TEGS system is available to the grid on demand. The Power availability factor (PAF) was computed to describe power availability. The PAF is computed as the fraction of time that the PV + storage system can deliver a given power to the grid. 
Moreover, PAF allows the examination of the availability of the combined power plant, as it measures how often a power plant supplies the rated power to the grid. A power plant with a 100\% PAF constantly supplies the grid with a given rated power. In the case of solar PV, the electricity generation suddenly drops when the sun no longer shines (owing to cloud cover or when the sun sets). This sudden drop in electricity generation reduces the PAF because the PV plant no longer generates at the rated power. Here, storage units can be used to charge whenever there is a low net demand for power in the grid and to discharge when there is a higher demand for power. Fig.~\ref{fig:TEGS_illustration} illustrates an example of how the storage unit can be used to shape the output to provide constant baseload power. Once the PV plant starts generating electricity over derated power (e.g., 20 MW), excess electricity is used to charge the storage unit. When the solar plant generates less electricity than the derated power, the storage unit starts discharging to satisfy the demand for electricity. 

Fig.~\ref{fig:TEGS_illustration} illustrates how the combined PV and storage system can provide constant baseload power, disregarding the electricity demand in the grid. However, when the system is connected to an electricity grid, it becomes significantly more complex. The system should not provide constant baseload power to the grid but should be able to supply the requested power to the grid system operator. 
In this study, the periods where the hypothetical PV + storage system cannot provide the requested electricity to the grid on demand are detected. These unwanted periods arise for the following reasons: 1) The PV system does not deliver the required power, and 2) The TEGS system cannot discharge the requested power as the State of Charge (SOC) is already zero. Such critical periods can be calculated as:
\begin{equation}
    \textrm{Percentage Not Available} = \frac{\textrm{When} ((PV_{gen} + TEGS_{discharge} < \textrm{Derated power}) + (SOC = 0))}{\textrm{8760 hours}},
\end{equation}
The system Percentage Not Available (PNA) gives information about the percentage of time during a full year the system cannot provide the requested electricity to the grid. On contrary, the PAF will be calculated as:
\begin{equation}
    PAF= 1 - PNA
\end{equation}

\section{Results and discussions}
\subsection{Power availability from PV and energy storage}
\begin{figure}[!ht]
    \centering
    \includegraphics[scale=1]{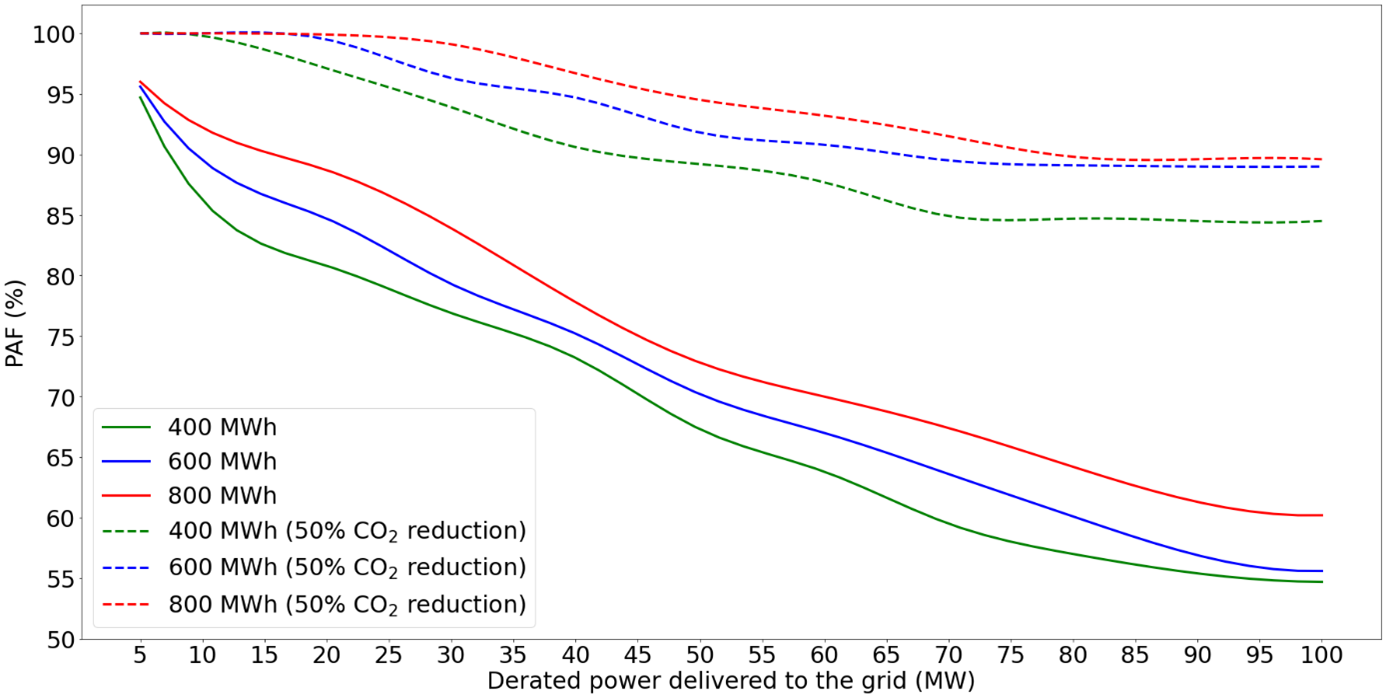}
    \caption{Example illustration where storage is used to provide a constant baseload power to the grid by charging when there is excess electricity and discharging when the PV plant does not generate electricity. The figure also illustrates the large variability in the daily generation profile of solar PV systems. }
    \label{fig:PAF_main}
\end{figure}

In Fig.~\ref{fig:PAF_main}, the yearly PAF from the hypothetical PV + TEGS system as a function of different discharge powers is provided. The more the system is derated, the more often it can deliver the required power. This is reasonable because the more the system is derated, the more the storage system can discharge at a lower-rated output for a longer period. If the discharge capacity of the storage system has a discharge capacity of 5 MW, the PV+TEGS system can deliver 5 MW to the grid approximately 95\% of the hours during the year. If the discharge capacity is 100 MW, the PV + TEGS system can deliver 100 MW to the grid approximately 55-60\% of the hours during the year. The percentage of time the system cannot deliver the required power also changes with different storage sizes, and there are fewer times during the year that the system cannot deliver the required power if the storage size is larger. 
Interestingly, there is a large difference between the scenarios with and without CO$_2$ reduction constraints. The electric power system dynamics change completely once CO$_2$ emissions were reduced by 50\%. In this case, the increased retirement of the NG makes the grid more dependent on the hypothetical PV + TEGS system, which results in the system supplying the necessary power 100\% of the time for a derated power between 5 MW and 20 MW for the TEGS unit with a storage capacity of 600 MWh and 800 MWh. This is remarkably higher than the baseline case, where the system cannot deliver the required power 5-15\% of the time for such derated powers. When modeling the PAF for the PV + TEGS system at the GW scale (see the Supplementary material), similar results are obtained. The PAF increases when lowering the derated power in both the baseline and the CO$_2$ reduction scenarios. In addition, the PAF is significantly higher when modeling the CO$_2$ reduction scenario compared to the baseline scenario. 
Although there are large differences in the PAF for the scenarios, the modelled PV + storage system has a significantly higher PAF than the average Capacity Factor (CF) (i.e., how often a power plant generates electricity at rated power) for single PV plants, which is approximately 24\% \cite{NRELATB}. Thus, the results show that combining storage with PV plants increases the dispatchability of the VRE technology and provides more reliable power output to the grid.

%ADD PLOTS HERE!
%

%ADD THE PLOT ASEGUN IS ASKING FOR HERE
%STATE THE FACT ASEGUN TALKS ABOUT THAT GAS IS NOT THAT IMPORTANT ANYMORE

The large difference between the CO$_2$ reduction scenarios indicates that the electric power system dynamic changes significantly if the grid must be less dependent on NG. Now the hypothetical PV + storage system plays a more important role in the grid because it does not emit any CO$_2$, and as such, the cost-minimization schedule of the CEM optimizes the grid to ensure that the PV + storage system can deliver the requested power to the grid more often during the year. 

\begin{figure}[!ht]
    \centering
    \includegraphics[scale=1]{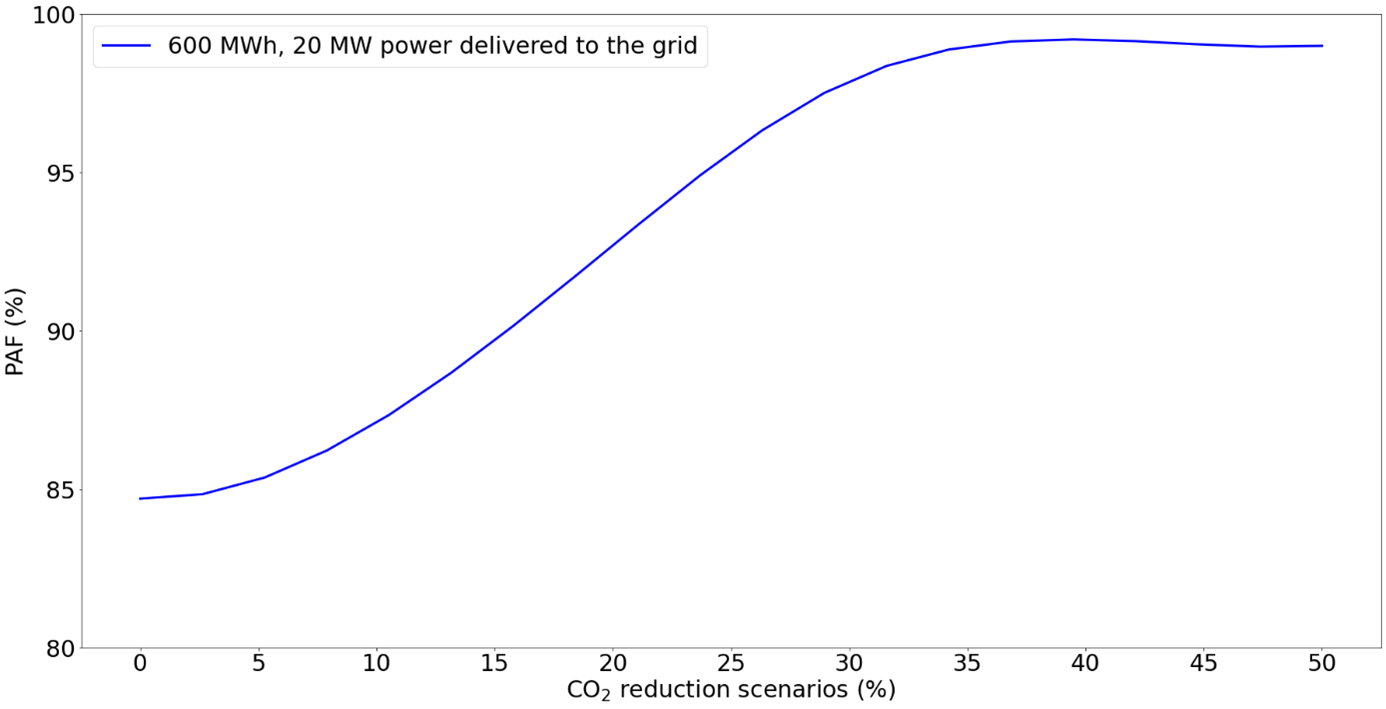}
    \caption{Percentage of time the hypothetical system can deliver the requested 20 MW power to the grid as a function of CO$_2$ reduction scenarios. The PAF increases as a function of CO$_2$ reduction as the power system must rely more on the hypothetical PV + storage system that does not emit any CO$_2$.}
    \label{fig:PAF_CO2red}
\end{figure}

Fig.~\ref{fig:PAF_CO2red} illustrates an example of how the PAF changes with different CO$_2$ reduction scenarios for a storage unit of 600 MWh that discharges 20 MW to the grid. Clearly, reducing CO$_2$ emissions results in a higher PAF. The maximum PAF was achieved at 50\% CO$_2$ reduction. Reducing CO$_2$ emissions requires the power system to be less dependent on fossil fuel technologies, such as NG, and thus must rely more on the power served by the PV + storage system. Modeling scenarios with more than 50\% CO$_2$ reduction results in an infeasible solution with the CEM optimization. Therefore, to further reduce the CO$_2$ emissions in the modeled grid, more renewable energy sources need to be installed in the current grid to replace NG sources.

\begin{figure}[!ht]
    \centering
    \includegraphics[scale=1]{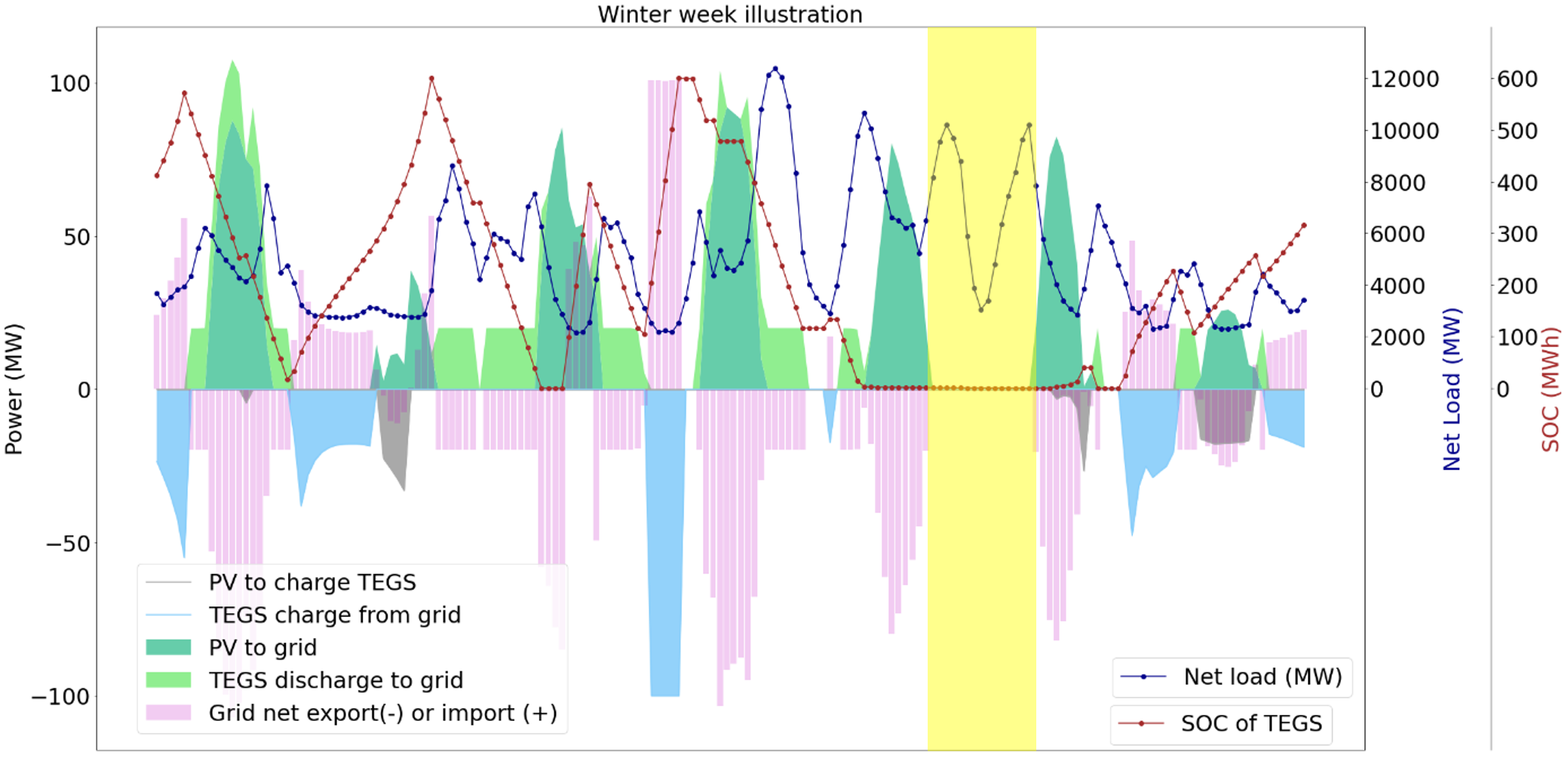}
    \includegraphics[scale=1]{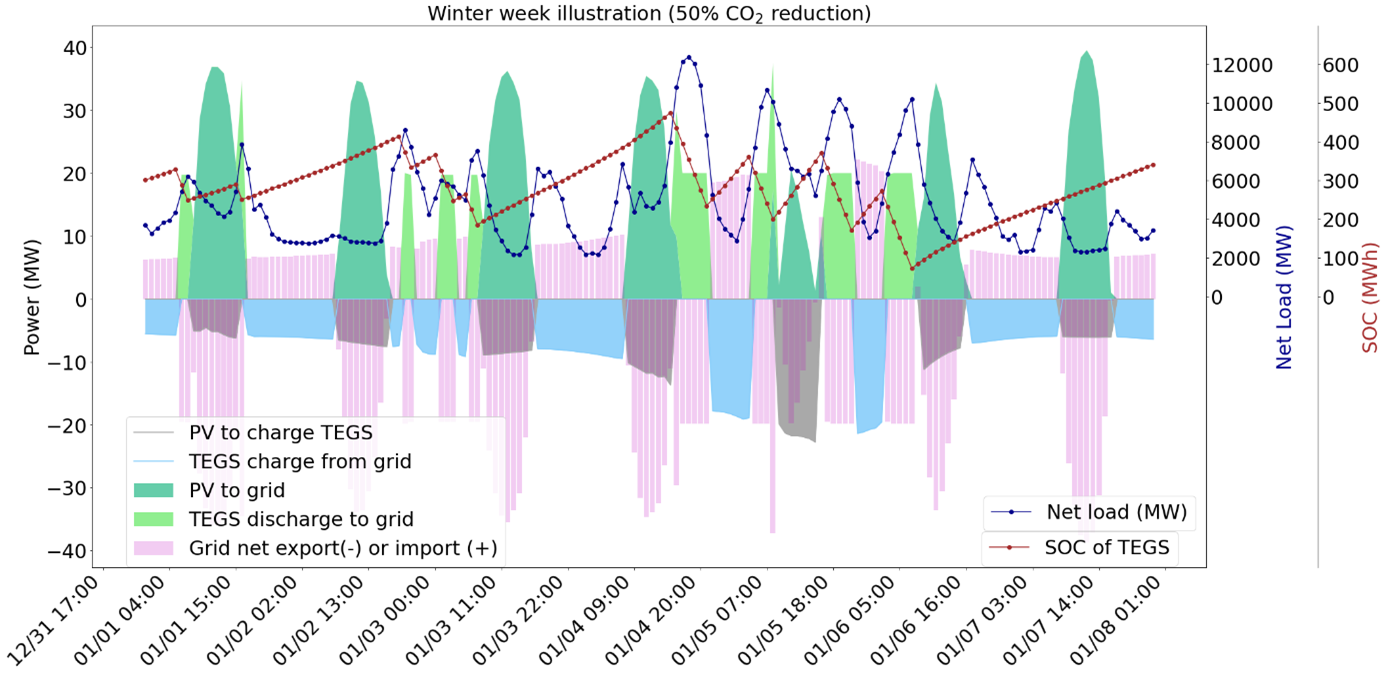}
    \caption{Winter week illustrations for a 100 MW PV plus a 600 MWh TEGS system that can discharge 20 MW. Baseline case without constraining CO$_2$ emissions (uppermost graph) and future scenario with reduced emissions lowermost. The baseline case shows a critical period during which the hypothetical system cannot deliver the requested power to the grid.  }
    \label{fig:Winter_eks}
\end{figure}

Fig.~\ref{fig:Winter_eks} shows how the hypothetical system operates in the different emissions scenarios (baseline and 50\% CO$_2$ reduction). Here, the hourly operation with a 600 MWh storage unit and 20 MW discharge capacity is shown during a typical winter week in January. The uppermost and lowermost figure show how the system operates in scenarios with and without constraining CO$_2$ emissions. 
From both graphs, TEGS is used frequently to discharge power to the grid whenever there is low electricity generation from the solar PV, which reduces the intermittency problem of PVs by increasing the number of hours the hypothetical system can deliver the required power to the grid. In addition, TEGS also charges power from the grid whenever there is a drop in the demand to increase the SOC, which illustrates the benefit of using storage units that are coupled with VRE resources. 
In the uppermost graph, a critical period is highlighted by the yellow squared region. Clearly, the grid has an increasing demand for electricity (blue line); however, the system cannot deliver the requested power because there is no generation from the PV system, and the storage system cannot discharge the required power to the grid because the SOC is already zero. However, considering the 50\% CO$_2$ reduction scenario, the system interacts differently with the grid, and it is clear that the CEM optimizes the hypothetical system to have more energy available at more times because the grid now is more dependent on the hypothetical system (because the grid can use less NG). For the particular example week, the PV + storage system can deliver the requested power to the grid at all times for the 50\% CO$_2$ reduction scenario. 

\begin{figure}[!ht]
    \centering
    \includegraphics[scale=1]{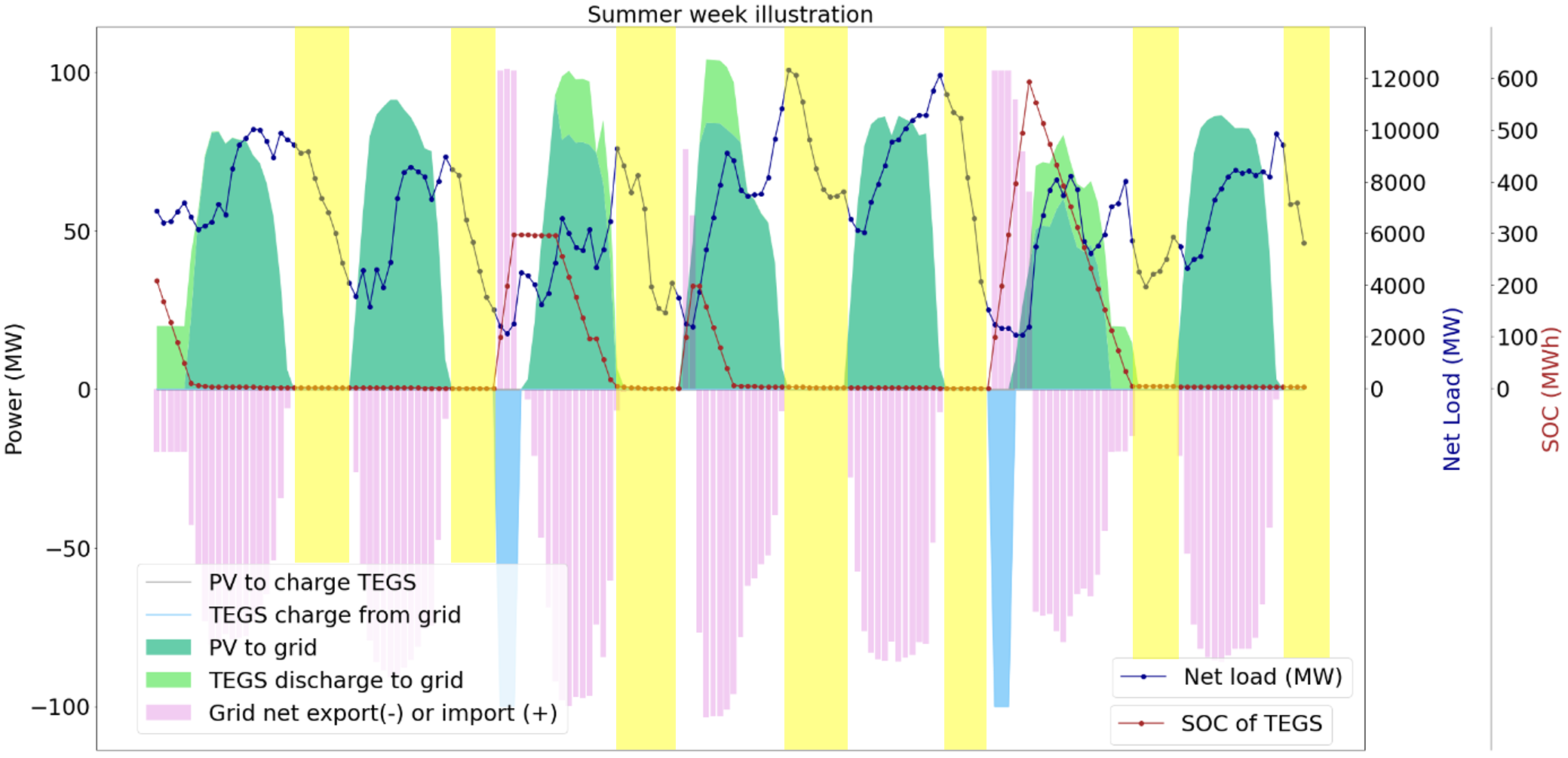}
    \includegraphics[scale=1]{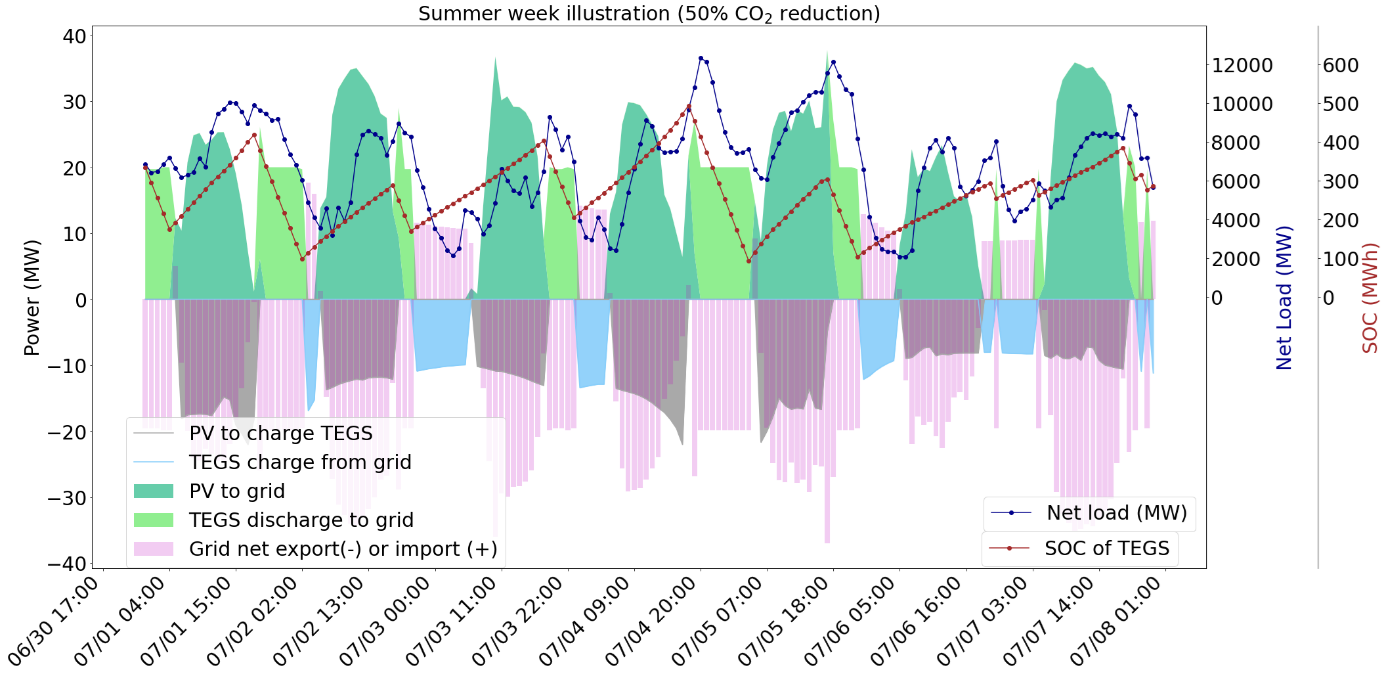}
    \caption{Illustrations of a typical summer week for a 100 MV PV plus 600 MWh TEGS system that can discharge 20 MW. Baseline case without constraining the CO$_2$ emissions (uppermost graph), and future scenario with reduced emissions lowermost. The PAF is significantly higher in the future scenarios. }
    \label{fig:Summer_eks}
    
\end{figure}

Like the winter week example shown in Fig.~\ref{fig:Winter_eks}, Fig.~\ref{fig:Summer_eks} shows how the hypothetical system operates with the grid during a typical summer week. Here, thanks to the higher solar availability, the system can deliver most of the required power to the grid using only the PV plant, and the storage system is used less frequently. However, for the baseline case, there are several periods in which the grid requests energy and the storage unit cannot provide sufficient power to the grid because the SOC is already zero. This is because the grid is mainly dependent on NG, which can supply power whenever there is low solar availability, and the TEGS system is only used to provide additional peak power when the demand suddenly increases. 
Under the 50\% CO$_2$ reduction scenario, the grid is much more dependent on the power from the storage unit whenever there is no PV generation. It is clear that instead of providing peak power to the grid, the PV system is used to charge the TEGS unit to a higher degree to ensure that the SOC is never zero and thus can discharge the derated power to the grid at all times when the PV system does not generate electricity. 

\section{Conclusion}

In this study, we analyzed how storage can be used to tackle the intermittency problem of VRE sources by increasing the dispatchability of a hypothetical PV plant. We modeled an existing electricity grid region in North America using a CEM and investigated how different storage configurations can reduce the number of periods in which a hypothetical PV + storage system cannot provide the required power to the grid. 

Because of the high capital cost of electrochemical batteries, a TES technology with a projected capital cost that fulfills the requirements (< US\$ 20 kWh$^{-1}$) to enable full decarbonization of the grid was considered. The power availability of the hypothetical system was modeled under different storage sizes and discharge capacities. Additionally, the optimization schedule was repeated under a hypothetical future scenario in which CO$_2$ emissions were constrained to be reduced by 50\%. In total, 66 different scenarios were modeled. To capture the high-resolution dependencies in the electricity generation balance, a full year with hourly resolution was optimized using the CEM. 

The results support the added value of using storage to increase the dispatchability of PV, as it significantly increases the PAF compared to PV systems alone. The percentage of time during the year the system could not deliver the required power to the grid decreased proportionally to the discharge capacity of the storage unit. In addition, increasing the storage size increases the energy availability, as more energy can be stored and thereafter discharged over a longer period when there is a demand for electricity in the grid. The findings were consistent when the PV plus TEGS system were evaluated at both the megawatt and gigawatt scale. 

Interestingly, there was a significant change in the electricity grid generation dynamics when the CO$_2$ emissions were reduced by 50\%. Here, because the grid can no longer rely on the same share of NG technology, the most cost-efficient grid is achieved when the PV + TEGS system is utilized to a higher degree, as these technologies do not emit any CO$_2$.
This shows that decreasing the maximum allowed GHG emissions in the grid significantly increases the value of using storage to increase the dispatchability of PV systems. 
We believe that our findings will provide increased knowledge to power system planners regarding how adding PV + storage systems to existing grids can contribute to the efficient stepwise decarbonization of power systems.

\subsection{Limitations and suggested future research}
This study presented an idealized representation of an existing grid (i.e., “brownfield” CEM approach) in the New England grid region. However, we are fully aware that the grid representation might not fully capture all details of the existing grid, and there can be differences between our abstract grid representation and the current real-world grid that is operated by ISO New England.  In addition to the transmission line between the existing grid and hypothetical PV + storage system, we modeled the current grid as a single-zone grid region without considering transmission losses or congestion between generators and demand. 
This study's major objectives are to evaluate a hypothetical PV + storage system's power availability and discuss the significance of integrating such technologies to assure a successful step-by-step decarbonization of the electric grid. Therefore, modeling the transmission lines between the existing generators is not considered, as it is outside the scope of this study and will significantly increase the computational intensity of the CEM.
The CEM is fully deterministic, assumes perfect foresight in planning and operational decisions, and does not account for uncertainty in VRE generation \cite{zhang2014review,salinas2020deepar,mashlakov2021assessing,eikeland2022probabilistic}. Therefore, this study does not aim to be used as a power planning tool for ISOs to assess the PAF of PV + storage systems in the day-ahead electricity market, but shows how storage, in general, will be a valuable technology to address the intermittency issue of VRE. A suggested future study is to frame the CEM to account for the uncertainty regarding the expected electricity generation from VRE sources. This will allow the use of the CEM as a decision-making tool for optimizing the management of the electricity grid in the day-ahead market.

\section{Acknowledgement}
O.F.E. and M.C. acknowledge the support from the research project “Transformation to a Renewable \& Smart Rural Power System Community (RENEW)”, connected to the Arctic Centre for Sustainable Energy (ARC) at UiT-the Arctic University of Norway through Grant No. 310026. We thank Maritsa Kissamitaki for designing Figure 1 and Figure 2. 

%Bibliography
\bibliographystyle{unsrt}  
\bibliography{templateArxiv}

\end{document}